\PassOptionsToPackage{dvipsnames}{xcolor} 
\documentclass[11pt,letter]{article}
\usepackage{amssymb,amsmath,amsthm,multirow,color,xcolor,cancel,bbm,tikz,mathrsfs}
\usepackage[textheight=24cm,textwidth=16cm]{geometry}
\tikzstyle{main node}=[draw,circle,inner sep=1,outer sep=2,thick,minimum size=12pt]

\newtheorem{lemma}{Lemma\setcounter{claimcounter}{0}}
\newtheorem{theorem}{Theorem\setcounter{claimcounter}{0}}

\newtheorem{question}{Question}

\newtheorem*{remark*}{Remark}

\newcounter{claimcounter}
\newtheorem*{claim*}{{\it Claim}}

\newcommand{\BF}[1]{{\boldmath{\bf #1}\unboldmath}}
\newcommand{\B}{\{0,1\}}

\title{Positive and negative cycles in Boolean networks\\{\it\normalsize In the memory of Ren\'e Thomas}}

\author{
Adrien Richard\footnote{Laboratoire I3S, CNRS, Universit\'e C\^ote d'Azur, France.
\tt{richard@unice.fr}}
}

\date{October 11, 2018}

\begin{document}

\maketitle

\begin{abstract}
We review and discuss some results about the influence of positive and negative feedback cycles in asynchronous Boolean networks. These results merge several ideas of Thomas: positive and negative feedback cycles have been largely emphasized by Thomas, through the so called Thomas' rules, and asynchronous Boolean networks have been introduced by Thomas as a model for the dynamics of gene networks, which is nowadays very popular. 
\end{abstract}

\section{Introduction}

We review and discuss some results, in discrete mathematics, that take their origin from two important contributions of Thomas: the {\em asynchronous logical description} of gene networks \cite{T73}, and the crucial influence of {\em positive and negative feedback cycles} in these networks~\cite{T81}.

\smallskip
Let us recall briefly the context. A gene network is usually roughly described, in a static way, by an {\em interaction graph}\footnote{which is sometime called {\em regulatory graph} or {\em influence graph}.}. This is a directed graph where arcs are signed positively or negatively. The vertices are the genes, and the presence of a positive (resp. negative) arc from one gene to another means that the product of the former has a positive (resp. negative) effect on the expression of the latter. In general, the interaction graph is rather well approximated while very little is known about the dynamics of the network ({\em i.e.}, the temporal evolution of the concentration of gene products). One is thus faced with the following difficult question: 
\[
\text{\em What can be said about the dynamics of a gene network from its interaction graph only?}
\]

\smallskip
To tackle this question from a mathematical point of view, one needs a model for the dynamics, and several models have been proposed. One of them is the asynchronous logical description introduced by Thomas and his coworkers \cite{T73,TA90,TK01}. It grasps the essential qualitative character of the dynamics while being consistent with the continuous descriptions based on non-linear ordinary differential equation systems \cite{S89,ST93}. In this setting, the time is {\em discrete} and concentrations evolve in {\em finite} domains. This allows the use of powerful computational tools to analysis the dynamics \cite{BCRG04,CNT12,A-J16,BL16}. Furthermore, the precision level of logical description is consistent with the lack of quantitative experimental data. For these reasons, this approach is more and more used nowadays \cite{A-J16}.  

\smallskip
Therefore, it makes sense to study the above question within the scope of the asynchronous logical description. This is what we will do here, by focusing on the crucial influence, emphasized by Thomas, of positive and negative feedback cycles of the interaction graph. To simplify the presentation, we will restrict to the Boolean case, firstly considered by Thomas. Gene products are then either present or absent. This case is sufficient for the phenomena we want to highlight. Many results remain valid for larger finite domains, and we will systematically discuss that. 

\section{Definitions}

\paragraph{Boolean networks.} 

In the Boolean case, the logical description consists in describing the dynamics from what is now called a Boolean network. Formally, a \BF{Boolean network} with $n$ components is a function 
\[
f:\B^n\to\B^n,\qquad x=(x_1,\dots,x_n)\mapsto f(x)=(f_1(x),\dots,f_n(x)). 
\]
The components $f_i$ of $f$ are thus Boolean functions, from $\B^n$ to $\B$. They are usually called the {\em local transition functions}. See Fig. \ref{fig:f} for an illustration. 

\begin{remark*}
In the context of gene networks, the interpretation is the following. There are $n$ genes, denoted from $1$ to $n$, and $\B^n$ is the set of possible {\em states} for the network: at state $x$, the product of gene $i$ is present if $x_i=1$ and absent if $x_i=0$. There are thus two possible concentration levels. The local transition function $f_i$ then gives the expression level of gene $i$ according to the state of the network: at state $x$, the gene $i$ is expressed if $f_i(x)=1$ and not expressed if $f_i(x)=0$.
\end{remark*}

\begin{figure}[h!]
\[
\begin{array}{c|c}
 x  & f(x)\\\hline
000 & 010 \\
001 & 000 \\
010 & 011 \\
011 & 000 \\
100 & 011 \\
101 & 010 \\
110 & 011 \\
111 & 111 
\end{array}
\qquad\qquad
\left\{
\begin{array}{l}
f_1(x)=x_1\land x_2\land x_3\\
f_2(x)=x_1\lor \overline{x_3}\\
f_3(x)=(x_2\land \overline{x_3})\lor (x_1\land \overline{x_2}\land\overline{x_3})\lor (x_1\land x_2\land x_3)
\end{array}
\right.
\]
\caption{\label{fig:f} A Boolean network $f$ with $3$ components, given in two equivalent ways: a table and a definition of the local transition functions with Boolean formula (in disjunctive normal form).} 
\end{figure}

\paragraph{Synchronous and asynchronous graphs.}

\begin{figure}[h!]
\[
\begin{array}{ccc}
\begin{array}{c}
\begin{array}{c}
\begin{tikzpicture}
\useasboundingbox (0,-0.5) rectangle (0,0.5);
\node (111) at (0,0){$111$};
\draw[->,thick] (111.-112) .. controls (-1,-1) and (1,-1) .. (111.-68);
\end{tikzpicture}
\end{array}
\begin{array}{c}
\begin{tikzpicture}
\node (000) at (150:1){$000$};
\node (010) at (30:1){$010$};
\node (011) at (-90:1){$011$};
\node (001) at (150:2.5){$001$};
\node (101) at (30:2.5){$101$};
\node (110) at (-110:2.1){$110$};
\node (100) at (-70:2.1){$100$};
\path[thick,->]
(000) edge[bend left=17] (010)
(010) edge[bend left=17] (011)
(011) edge[bend left=17] (000)
(001) edge (000)
(101) edge (010)
(110) edge (011)
(100) edge (011)
;
\end{tikzpicture}
\end{array}
\end{array}
&\qquad&
\begin{array}{c}
\begin{tikzpicture}
\node (000) at (-1.5,-1.5){$000$};
\node (001) at (-0.5,-0.5){$001$};
\node (010) at (-1.5, 0.5){$010$};
\node (011) at (-0.5, 1.5){$011$};
\node (100) at ( 0.5,-1.5){$100$};
\node (101) at ( 1.5,-0.5){$101$};
\node (110) at ( 0.5, 0.5){$110$};
\node (111) at ( 1.5, 1.5){$111$};
\path[thick,->,draw,black]
(000) edge (010)
(001) edge (000)
(010) edge[bend left=10] (011)
(011) edge (001)
(011) edge[bend left=10] (010)
(100) edge (000)
(100) edge (110)
(100) edge[bend left=10] (101)
(101) edge (001)
(101) edge (111)
(101) edge[bend left=10] (100)
(110) edge (010)
(110) edge (111)
;
\end{tikzpicture}
\end{array}
\\~\\
\Sigma(f)&&\Gamma(f)
\end{array}
\]
\caption{\label{fig:dynamics} The synchronous and asynchronous graphs of the Boolean network $f$ of Fig. \ref{fig:f}.}
\end{figure}
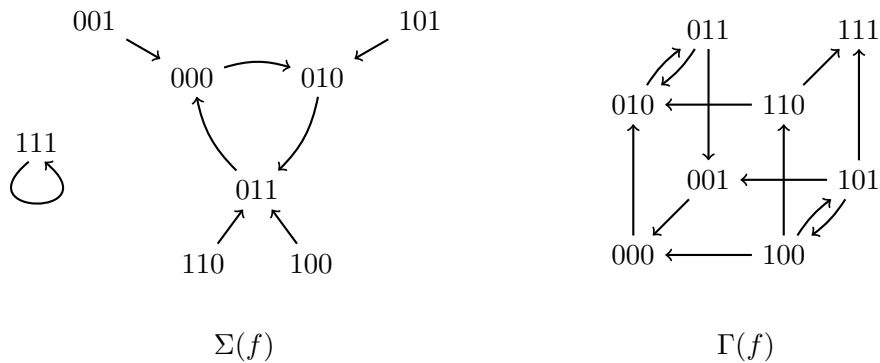

The \BF{synchronous graph} of $f$, denoted $\Sigma(f)$, is the directed graph with vertex set $\B^n$ and an arc  $x\to f(x)$ for each state $x$. It describes the {\em synchronous dynamics}: if $x^t$ is the state at time $t$, then $x^{t+1}=f(x^t)$. Thus, at each time step, many gene products may appear or disappear synchronously. As argued by Thomas,  this is not realistic and not consistent with continuous models, contrary to the {\em asynchronous dynamics} he introduced (see \cite{TK01} and the reference therein). Formally, the asynchronous dynamics is described by the \BF{asynchronous graph} $\Gamma(f)$ defined as follows: the vertex set is $\B^n$ and, for every state $x$ and every component $i$ such that $f_i(x)\neq x_i$, there is an arc (or transition) $x\to \bar x^i$, where $\bar x^i$ is the state obtained from $x$ by flipping the $i$th component. Since several transitions can start from the same state, this dynamics is {\em non-deterministic}: if $x^t$ is the state at time $t$, then $x^{t+1}$ may be {\em any} out-neighbor of $x^t$ in $\Gamma(f)$. See Fig. \ref{fig:dynamics} for an illustration. 

\begin{remark*}
In the context of gene networks, the interpretation is the following. If $f_i(x)\neq x_i$ then the concentration level of the product of gene $i$ is called to change: either the gene $i$ is not expressed and its product is present ($f_i(x)<x_i$), or the gene $i$ is expressed and its product is absent ($f_i(x)>x_i$). If $f_i(x)=x_i$, then the concentration level of the product of $i$ is stable (gene expressed and product present, or gene not expressed and product absent). Hence, each transition of $\Gamma(f)$ starting from $x$ corresponds to the evolution of one concentration level, which is called to change at state $x$. This is suited for a situation where the synchronous evolution of several concentration levels is unlikely. 
\end{remark*}

\paragraph{Local and global interaction graphs.}

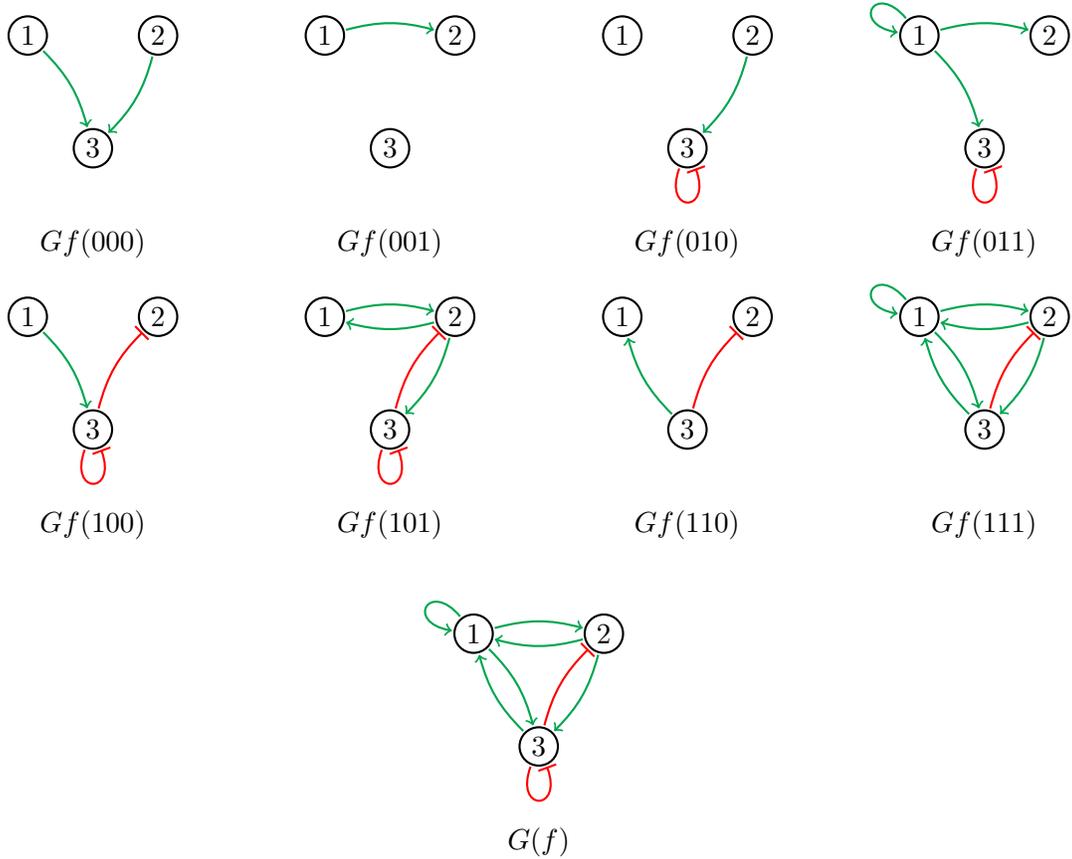
\begin{figure}[h!]
\[
\begin{array}{c}
\begin{array}{cccc}
\begin{tikzpicture}
\node[outer sep=1,inner sep=2,circle,draw,thick] (1) at (150:1){$1$};
\node[outer sep=1,inner sep=2,circle,draw,thick] (2) at (30:1){$2$};
\node[outer sep=1,inner sep=2,circle,draw,thick] (3) at (270:1){$3$};
\draw[white,->,thick] (1.128) .. controls (140:1.9) and (160:1.9) .. (1.172);
\draw[white,->,thick] (2.8) .. controls (20:1.9) and (40:1.9) .. (2.52);
\draw[white,-|,thick] (3.-112) .. controls (-100:1.9) and (-80:1.9) .. (3.-68);
\path[->,thick]
(1) edge[white,bend left=15] (2)
(2) edge[Green,bend left=15] (3)
(3) edge[white,bend left=15] (1)
(1) edge[Green,bend left=15] (3)
(3) edge[-|,white,bend left=15] (2)
(2) edge[white,bend left=15] (1)
;
\end{tikzpicture}
&
\begin{tikzpicture}
\node[outer sep=1,inner sep=2,circle,draw,thick] (1) at (150:1){$1$};
\node[outer sep=1,inner sep=2,circle,draw,thick] (2) at (30:1){$2$};
\node[outer sep=1,inner sep=2,circle,draw,thick] (3) at (270:1){$3$};
\draw[white,->,thick] (1.128) .. controls (140:1.9) and (160:1.9) .. (1.172);
\draw[white,->,thick] (2.8) .. controls (20:1.9) and (40:1.9) .. (2.52);
\draw[white,-|,thick] (3.-112) .. controls (-100:1.9) and (-80:1.9) .. (3.-68);
\path[->,thick]
(1) edge[Green,bend left=15] (2)
(2) edge[white,bend left=15] (3)
(3) edge[white,bend left=15] (1)
(1) edge[white,bend left=15] (3)
(3) edge[-|,white,bend left=15] (2)
(2) edge[white,bend left=15] (1)
;
\end{tikzpicture}
&
\begin{tikzpicture}
\node[outer sep=1,inner sep=2,circle,draw,thick] (1) at (150:1){$1$};
\node[outer sep=1,inner sep=2,circle,draw,thick] (2) at (30:1){$2$};
\node[outer sep=1,inner sep=2,circle,draw,thick] (3) at (270:1){$3$};
\draw[white,->,thick] (1.128) .. controls (140:1.9) and (160:1.9) .. (1.172);
\draw[white,->,thick] (2.8) .. controls (20:1.9) and (40:1.9) .. (2.52);
\draw[red,-|,thick] (3.-112) .. controls (-100:1.9) and (-80:1.9) .. (3.-68);
\path[->,thick]
(1) edge[white,bend left=15] (2)
(2) edge[Green,bend left=15] (3)
(3) edge[white,bend left=15] (1)
(1) edge[white,bend left=15] (3)
(3) edge[-|,white,bend left=15] (2)
(2) edge[white,bend left=15] (1)
;
\end{tikzpicture}
&
\begin{tikzpicture}
\node[outer sep=1,inner sep=2,circle,draw,thick] (1) at (150:1){$1$};
\node[outer sep=1,inner sep=2,circle,draw,thick] (2) at (30:1){$2$};
\node[outer sep=1,inner sep=2,circle,draw,thick] (3) at (270:1){$3$};
\draw[Green,->,thick] (1.128) .. controls (140:1.9) and (160:1.9) .. (1.172);
\draw[white,->,thick] (2.8) .. controls (20:1.9) and (40:1.9) .. (2.52);
\draw[red,-|,thick] (3.-112) .. controls (-100:1.9) and (-80:1.9) .. (3.-68);
\path[->,thick]
(1) edge[Green,bend left=15] (2)
(2) edge[white,bend left=15] (3)
(3) edge[white,bend left=15] (1)
(1) edge[Green,bend left=15] (3)
(3) edge[-|,white,bend left=15] (2)
(2) edge[white,bend left=15] (1)
;
\end{tikzpicture}
\\
Gf(000)&Gf(001)&Gf(010)&Gf(011)
\end{array}
\\
\begin{array}{cccc}
\begin{tikzpicture}
\node[outer sep=1,inner sep=2,circle,draw,thick] (1) at (150:1){$1$};
\node[outer sep=1,inner sep=2,circle,draw,thick] (2) at (30:1){$2$};
\node[outer sep=1,inner sep=2,circle,draw,thick] (3) at (270:1){$3$};
\draw[white,->,thick] (1.128) .. controls (140:1.9) and (160:1.9) .. (1.172);
\draw[white,->,thick] (2.8) .. controls (20:1.9) and (40:1.9) .. (2.52);
\draw[red,-|,thick] (3.-112) .. controls (-100:1.9) and (-80:1.9) .. (3.-68);
\path[->,thick]
(1) edge[white,bend left=15] (2)
(2) edge[white,bend left=15] (3)
(3) edge[white,bend left=15] (1)
(1) edge[Green,bend left=15] (3)
(3) edge[-|,red,bend left=15] (2)
(2) edge[white,bend left=15] (1)
;
\end{tikzpicture}
&
\begin{tikzpicture}
\node[outer sep=1,inner sep=2,circle,draw,thick] (1) at (150:1){$1$};
\node[outer sep=1,inner sep=2,circle,draw,thick] (2) at (30:1){$2$};
\node[outer sep=1,inner sep=2,circle,draw,thick] (3) at (270:1){$3$};
\draw[white,->,thick] (1.128) .. controls (140:1.9) and (160:1.9) .. (1.172);
\draw[white,->,thick] (2.8) .. controls (20:1.9) and (40:1.9) .. (2.52);
\draw[red,-|,thick] (3.-112) .. controls (-100:1.9) and (-80:1.9) .. (3.-68);
\path[->,thick]
(1) edge[Green,bend left=15] (2)
(2) edge[Green,bend left=15] (3)
(3) edge[white,bend left=15] (1)
(1) edge[white,bend left=15] (3)
(3) edge[-|,red,bend left=15] (2)
(2) edge[Green,bend left=15] (1)
;
\end{tikzpicture}
&
\begin{tikzpicture}
\node[outer sep=1,inner sep=2,circle,draw,thick] (1) at (150:1){$1$};
\node[outer sep=1,inner sep=2,circle,draw,thick] (2) at (30:1){$2$};
\node[outer sep=1,inner sep=2,circle,draw,thick] (3) at (270:1){$3$};
\draw[white,->,thick] (1.128) .. controls (140:1.9) and (160:1.9) .. (1.172);
\draw[white,->,thick] (2.8) .. controls (20:1.9) and (40:1.9) .. (2.52);
\draw[white,-|,thick] (3.-112) .. controls (-100:1.9) and (-80:1.9) .. (3.-68);
\path[->,thick]
(1) edge[white,bend left=15] (2)
(2) edge[white,bend left=15] (3)
(3) edge[Green,bend left=15] (1)
(1) edge[white,bend left=15] (3)
(3) edge[-|,red,bend left=15] (2)
(2) edge[white,bend left=15] (1)
;
\end{tikzpicture}
&
\begin{tikzpicture}
\node[outer sep=1,inner sep=2,circle,draw,thick] (1) at (150:1){$1$};
\node[outer sep=1,inner sep=2,circle,draw,thick] (2) at (30:1){$2$};
\node[outer sep=1,inner sep=2,circle,draw,thick] (3) at (270:1){$3$};
\draw[Green,->,thick] (1.128) .. controls (140:1.9) and (160:1.9) .. (1.172);
\draw[white,->,thick] (2.8) .. controls (20:1.9) and (40:1.9) .. (2.52);
\draw[white,-|,thick] (3.-112) .. controls (-100:1.9) and (-80:1.9) .. (3.-68);
\path[->,thick]
(1) edge[Green,bend left=15] (2)
(2) edge[Green,bend left=15] (3)
(3) edge[Green,bend left=15] (1)
(1) edge[Green,bend left=15] (3)
(3) edge[-|,red,bend left=15] (2)
(2) edge[Green,bend left=15] (1)
;
\end{tikzpicture}
\\
Gf(100)&Gf(101)&Gf(110)&Gf(111)
\end{array}
\\~\\
\begin{array}{c}
\begin{tikzpicture}
\node[outer sep=1,inner sep=2,circle,draw,thick] (1) at (150:1){$1$};
\node[outer sep=1,inner sep=2,circle,draw,thick] (2) at (30:1){$2$};
\node[outer sep=1,inner sep=2,circle,draw,thick] (3) at (270:1){$3$};
\draw[Green,->,thick] (1.128) .. controls (140:1.9) and (160:1.9) .. (1.172);
\draw[white,->,thick] (2.8) .. controls (20:1.9) and (40:1.9) .. (2.52);
\draw[red,-|,thick] (3.-112) .. controls (-100:1.9) and (-80:1.9) .. (3.-68);
\path[->,thick]
(1) edge[Green,bend left=15] (2)
(2) edge[Green,bend left=15] (3)
(3) edge[Green,bend left=15] (1)
(1) edge[Green,bend left=15] (3)
(3) edge[-|,red,bend left=15] (2)
(2) edge[Green,bend left=15] (1)
;
\end{tikzpicture}
\\
G(f)
\end{array}
\end{array}
\]
\caption{\label{fig:graphs} The local interaction graphs and the global interaction graph of the Boolean network $f$ of Fig. \ref{fig:f}. Green arrows correspond to positive arcs, and T-end red arrows correspond to negative arcs. This convention is used throughout the paper.}
\end{figure}

We will now define the {\em interaction graph} of $f$. For that we need some useful preliminary definitions. The \BF{discrete partial derivative} of $f_i$ with respect to the component $j$ is the function $f_{ij}$ from $\B^n$ to $\{-1,0,1\}$ defined by:
\[
f_{ij}(x):=f_i(x_1,\dots,x_{j-1},1,x_{j+1},\dots,x_n)-f_i(x_1,\dots,x_{j-1},0,x_{j+1},\dots,x_n).
\] 
Hence, $f_{ij}$ is the discrete analogue of $\partial f_i/\partial x_j$. The \BF{local interaction graph} of $f$ at state $x$, denoted $Gf(x)$, is then the signed directed graph defined as follows: the vertex set is $[n]:=\{1,\dots,n\}$ and, for all $i,j\in [n]$ (not necessarily distinct), there is a positive (resp. negative) arc $j\to i$ if $f_{ij}(x)$ is positive (resp. negative). Finally, the \BF{global interaction graph} of $f$, denoted $G(f)$, is the union of all the local interaction graphs: the vertex set is $[n]$ and, for all $i,j\in [n]$ (not necessarily distinct), there is a positive (resp. negative) arc $j\to i$ if $f_{ij}(x)$ is positive (resp. negative) for {\em at least one} state $x$. Note that, contrary to the local interactions graphs, the global interaction graph $G(f)$ can have both a positive and a negative arc from one vertex to another. In that case, the sign of the interaction depends on the state of the network. In the following, $G(f)$ is simply called the \BF{interaction graph} of $f$. In a signed directed graph, the sign of a cycle is the product of the sign of its arcs. Thus a cycle is positive if and only if it contains an even number of negative arcs\footnote{In a signed or unsigned directed graph, paths and cycles are always directed and without repeated vertices, and in the signed case, they do not contain both a positive and a negative arc from one vertex to another.}. See Fig. \ref{fig:graphs} for an illustration. 

\begin{remark*}
In the context of gene networks, the interpretation is the following. If $f_{ij}(x)>0$ then it means that, at state $x$, an increase of the concentration level of the product of gene $j$ switches on the expression level of gene $i$, and thus, at state $x$, the product of $j$ has a positive effect on the expression of $i$, which is reported by a positive arc $j\to i$ in the local interaction graph $Gf(x)$. Symmetrically, $f_{ij}(x)<0$ then it means that, at state $x$, an increase of the concentration level of the product of gene $j$ switches off the expression level of gene $i$, and thus, at state $x$, the product of $j$ has a negative effect on the expression of $i$, which is reported by a negative arc $j\to i$ in the local interaction graph $G(f,x)$. The interaction graph $G(f)$ then simply collects all the interactions that can be seen locally.  
\end{remark*}

\paragraph{Attractors and fixed points.}

At this stage, the question of the introduction can be restated more precisely as follows: 
\[
\text{\em What can be said about $\Gamma(f)$ according to $G(f)$?} 
\]
This question is still very large, and we will thus focus on special features of $\Gamma(f)$, essentially the \BF{attractors} of $\Gamma(f)$. These are defined as the smallest subsets $X\subseteq \B^n$ (with respect to the inclusion relation) such that there is no transition $x\to y$ in $\Gamma(f)$ with $x\in X$ and $y\not\in X$. Thus attractors are the smallest subsets of states that we cannot leave in the asynchronous dynamics. Equivalently, an attractor is the vertex set of a {\em strongly connected component} of $\Gamma(f)$ without leaving transition. Thus, inside an attractor $X$ of size $|X|\geq 2$, the dynamics necessarily describes sustained oscillations (that are possibly aperiodic). Attractors of size at least two are then called \BF{cyclic attractors}. Furthermore, $X=\{x\}$ is an attractor of size one if and only if $x$ is a \BF{fixed point} of $f$, that is, $x=f(x)$. Thus attractors of size one and fixed points can be identified and correspond to \BF{stable states}.  
%
%
For instance, the asynchronous graph in Fig.~\ref{fig:dynamics} contains two attractors: the fixed point 111 and the cyclic attractor $\{000,010,001,011\}$. Many results we will discuss concern the following more precise question:
\[
\text{\em What can be said about the number of fixed points of $f$ according to $G(f)$?} 
\]

\begin{remark*}
In the context of gene networks, fixed points describe stable patterns of gene expression and are often related to particular cellular processes. In particular, the presence of multiple fixed points should accounts for differentiation processes \cite{TK01a,TK01}. Therefore, the number of fixed points is an important dynamical property. The above question is thus of special interest. 
\end{remark*}

\paragraph{Paths and geodesic paths.}

We will also be interested in reachability properties. In a directed graph, the length of a path is the number of arcs it contains. The \BF{Hamming distance} $d(x,y)$ between two states $x,y\in\B^n$ is the number of components $i\in [n]$ such that $x_i\neq y_i$. In the asynchronous graph $\Gamma(f)$, a path from $x$ to $y$ is necessarily of length at least $d(x,y)$. We then say that a path of $\Gamma(f)$ from $x$ to $y$ is a \BF{geodesic} if its length is exactly $d(x,y)$. Equivalently, a geodesic is a path along which each component changes at most one time.  

\begin{remark*}
All the notions introduced above have natural extensions to the {\em discrete} case, where $\B^n$ is replaced by the product of $n$ finite integer intervals, see \cite{RC07} for instance. In the following, all the graphs we consider are directed, and we thus omit this precision. We also assume some basic notions about graphs ({\em e.g.} subgraph, strongly connected component), following the terminology of \cite{BG08}.
\end{remark*}

\section{The importance of feedback cycles}

So, what can be said about the dynamics described by $f$ (synchronous or aysnchronous) according to $G(f)$? To start, it is natural to make strong assumptions on $G(f)$ and see what happens. For instance, what can be said when $G(f)$ is {\em acyclic}, that is, does not contain any cycle? A rather complete answer has been obtained by Robert \cite{R80,R86,R95} in the discrete setting. In the Boolean case the statement is the following. 

\begin{theorem}[\BF{Robert}]\label{thm:robert}
Let $f:\B^n\to\B^n$ and suppose that $G(f)$ is acyclic. Then
\begin{enumerate}
\item $f$ has a unique fixed point, say $x$.
\item $\Sigma(f)$ has a path of length at most $n$ from every state $y$ to $x$. 
\item $\Gamma(f)$ is acyclic and has a geodesic path from every state $y$ to $x$.
\end{enumerate}
\end{theorem}

Thus, in the acyclic case, there is both a synchronous and asynchronous convergence toward a unique fixed point: both dynamics are very simple. This means that complex behaviors need the presence of feedback cycles: {\em ``feedback cycles are the engines of complexity''}. Note that point~{\em 3.} implies that $\Gamma(f)$ has a unique attractor, which is the unique fixed point~of~$f$. 

\smallskip
Robert's theorem is not the end of the story concerning the acyclic case. Shih and Dong \cite{SD05} widely generalized the first point, as follows.

\begin{theorem}[\BF{Shih and Dong}]
Let $f:\B^n\to\B^n$ and suppose that $Gf(x)$ is acyclic for every $x\in\{0,1\}^n$. Then $f$ has a unique fixed point.
\end{theorem}

To see why this theorem generalizes the first point in Robert's theorem, it is sufficient to remark that if $G(f)$ is acyclic, then so is $Gf(x)$ for every $x$, since each $Gf(x)$ is a subgraph of $G(f)$. Hence, Shih and Dong obtained the same conclusion (a unique fixed point) under a weaker assumption. However, under this weaker assumption, points {\em 2.} and {\em 3.} are lost. The $4$-component example in \cite{SD05} shows that.  Note that Shih and Dong's theorem has been conjectured and presented as a Boolean analogue of the Jacobian conjecture in algebraic geometry in \cite{SH99}.   

\smallskip
Shih and Dong's theorem has been extended to the discrete case in \cite{R08} and a generalization in the Boolean case is given in \cite{R15}. Roughly speaking, it shows that $f$ still has a unique fixed point when there are ``few'' states $x$ such that $Gf(x)$ has a ``short'' cycle. 

\begin{theorem}\label{thm:shihdong}
Let $f:\B^n\to\B^n$ and suppose that, for every $1\leq \ell\leq n$, there is less than $2^\ell$ states $x\in\B^n$ such that $Gf(x)$ has a cycle of length $\leq \ell$. Then $f$ has a unique fixed~point.
\end{theorem}

By Robert's theorem, when $G(f)$ is acyclic, there exists, from every initial state, a short asynchronous path toward the unique fixed point of $f$. But this does not prevent from the existence of long asynchronous paths. Indeed, Domshlak \cite{D02} proved (in a different context) that it may exist two states $x$ and $y$ such that $\Gamma(f)$ has paths from $x$ to $y$,  all of exponential lengths with respect to $n$. 

\begin{theorem}
For every $n\geq 8$ there exists $f:\B^n\to\B^n$ with $G(f)$ acyclic, and two states $x,y\in\B^n$, such that $\Gamma(f)$ has paths from $x$ to $y$, all of length at least $1.5^{\frac{n}{2}}$.
\end{theorem}
 
\section{Thomas' rules}

Robert's theorem shows that cycles are necessary for complex behaviors, and Thomas emphasized the fact that it really matters to consider two types of cycles: the positive and negative ones. These two types of cycles have rather opposite dynamical influences. 

\begin{figure}[h!]
\[
\begin{array}{cc}
\begin{array}{c|c}
x & f(x)\\\hline
00 & 11\\
01 & 01\\
10 & 10\\
11 & 11
\end{array}
\quad&
\begin{array}{ccccc}
\begin{array}{c}
\begin{tikzpicture}
\node[outer sep=1,inner sep=2,circle,draw,thick] (1) at (0,0){$1$};
\node[outer sep=1,inner sep=2,circle,draw,thick] (2) at (1.5,0){$2$};
\path[red,-|,thick]
(1) edge[bend left=35] (2)
(2) edge[bend left=35] (1)
;
\end{tikzpicture}
\end{array}
&\quad&
\begin{array}{c}
\begin{tikzpicture}
\useasboundingbox (-0.3,-0.3) rectangle (1.6,1.6);
\node (00) at (0,0){$00$};
\node (01) at (0,1.3){\BF{$01$}};
\node (10) at (1.3,0){\BF{$10$}};
\node (11) at (1.3,1.3) {$11$};
\path[thick,->]
(00) edge (01)
(00) edge (10)
(11) edge (01)
(11) edge (10)
;
\end{tikzpicture}
\end{array}
&\quad&
\begin{array}{c}
\begin{tikzpicture}
\useasboundingbox (-0.3,-0.3) rectangle (1.6,1.6);
\node (00) at (0,0){$00$};
\node (01) at (0,1.3){$01$};
\node (10) at (1.3,0){$10$};
\node (11) at (1.3,1.3) {$11$};
\draw[->,thick] (01) .. controls (-1,1.3) and (0,2.3) .. (01);
\draw[->,thick] (10) .. controls (2.3,0) and (1.3,-1) .. (10);
\path[thick,->]
(00) edge[bend left=15] (11)
(11) edge[bend left=15] (00)
;
\end{tikzpicture}
\end{array}
\\[10mm]
G(f) &&\Gamma(f)&&\Sigma(f)
\end{array}
\end{array}
\]
\caption{\label{fig:positive_cycle} A two-component Boolean network whose interaction graph is a positive cycle.}
\end{figure}
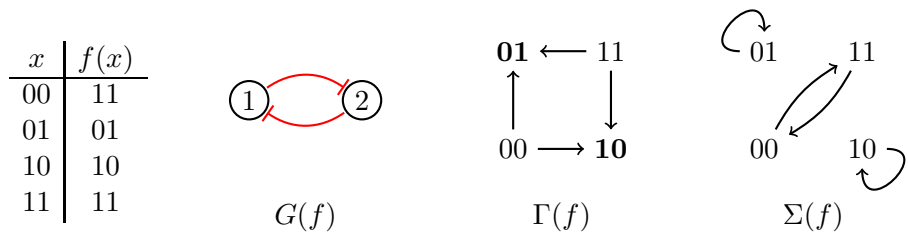

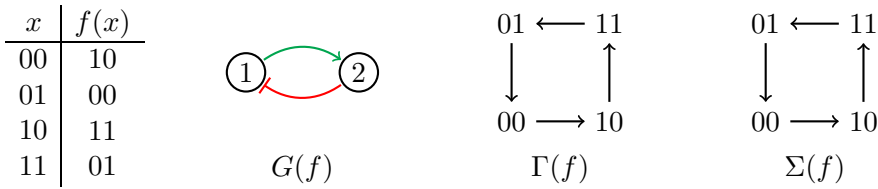
\begin{figure}[h!]
\[
\begin{array}{cc}
\begin{array}{c|c}
x & f(x)\\\hline
00 & 10\\
01 & 00\\
10 & 11\\
11 & 01
\end{array}
\quad
&
\begin{array}{ccccc}
\begin{array}{c}
\begin{tikzpicture}
\node[outer sep=1,inner sep=2,circle,draw,thick] (1) at (0,0){$1$};
\node[outer sep=1,inner sep=2,circle,draw,thick] (2) at (1.5,0){$2$};
\path[thick]
(1) edge[Green,->,bend left=35] (2)
(2) edge[red,-|,bend left=35] (1)
;
\end{tikzpicture}
\end{array}
&\quad&
\begin{array}{c}
\begin{tikzpicture}
\node (00) at (0,0){$00$};
\node (01) at (0,1.3){$01$};
\node (10) at (1.3,0){$10$};
\node (11) at (1.3,1.3) {$11$};
\path[thick,->]
(00) edge (10)
(10) edge (11)
(11) edge (01)
(01) edge (00)
;
\end{tikzpicture}
\end{array}
&\quad&
\begin{array}{c}
\begin{tikzpicture}
\node (00) at (0,0){$00$};
\node (01) at (0,1.3){$01$};
\node (10) at (1.3,0){$10$};
\node (11) at (1.3,1.3) {$11$};
\path[thick,->]
(00) edge (10)
(10) edge (11)
(11) edge (01)
(01) edge (00)
;
\end{tikzpicture}
\end{array}
\\
G(f) &&\Gamma(f)&&\Sigma(f)
\end{array}
\end{array}
\]
\caption{\label{fig:negative_cycle} A two-component Boolean network whose interaction graph is a negative cycle.}
\end{figure}

\smallskip
To make a first distinction, suppose that $G(f)$ is itself a single isolated cycle. If $G(f)$ is positive, then $\Gamma(f)$ has exactly two attractors, which are fixed points, while if $G(f)$ is negative, then $\Gamma(f)$ has a unique cyclic attractor, of size $2n$, that describes periodic and deterministic sustained oscillations (each state of the attractor has a unique out-neighbor). See Fig.~\ref{fig:positive_cycle} and Fig.~\ref{fig:negative_cycle} for an illustration. By Robert's theorem, we can conclude that positive cycles are the simplest interaction graphs allowing the presence of multiple stable states, while negative cycles are the simplest interaction graphs allowing the presence of sustained oscillations (cyclic attractors). A comprehensive description of the dynamics of isolated cycles is given in \cite{RMCT03} for the asynchronous case, and in \cite{DNS12} for the synchronous case (the case of two intersecting cycles is also well studied, see \cite{DENST11,DENS11,RMCT03,MNRS13}). 

\smallskip
However, the presence of a positive (resp. negative) cycle in $G(f)$ is {\em not} a sufficient condition for the presence of multiple fixed points (resp. sustained oscillations). For instance, if $G(f)$ consists of a positive cycle and a negative cycle that intersect each other at a unique vertex, then $\Gamma(f)$ has a unique attractor, which is a fixed point \cite{MNRS13}. 

\smallskip
These observations leaded Thomas \cite{T81} to make two general conjectures, sometime called rules since statements are not given in a specific mathematical framework, the following:
\begin{enumerate}
\item
First rule: the presence of a positive cycle in the interaction graph of a dynamical system is a necessary condition for multiple stable steady states. 
\item
Second rule: the presence of a negative cycle in the interaction graph of a dynamical system is a necessary condition for permanent periodic behaviors. 
\end{enumerate}

\smallskip
These two rules gave rise to several mathematical theorems, both in the continuous and discrete settings. Concerning the continuous case, we refer the reader to the paper of Kaufman and Soul\'e in this Special Issue (see also \cite{CD02,S03,S06} and the references therein). The present paper discusses the Boolean case. The interpretation of Thomas' rules is then very natural: the dynamical system is $\Gamma(f)$, the interaction graph is $G(f)$, the stable steady states are the fixed points of $f$, and the permanent periodic behaviors are produced by the cyclic attractors of $\Gamma(f)$. 

\section{Thomas' first rule and beyond}

Aracena \cite{A08} (see also \cite{ADG04a}) proved the following natural Boolean version of Thomas' first rule. 

\begin{theorem}[\BF{Thomas' first rule, global version}]\label{thm:thomas1}
Let $f:\B^n\to\B^n$ and suppose that $G(f)$ has no positive cycle. Then $f$ has at most one fixed point. 
\end{theorem}

Aracena also proved, in the same paper, the following theorem, with a stronger conclusion under a stronger assumption (it is an easy exercise to deduce from it the above theorem, by considering the strongly connected components in the topological order).

\begin{theorem}
Let $f:\B^n\to\B^n$ and suppose that $G(f)$ is strongly connected, has at least one arc, and has no positive cycle. Then $f$ has no fixed point. 
\end{theorem}

Theorem~\ref{thm:thomas1} has many extensions that we discuss below. 
Actually, these extensions are obtained with the following slightly stronger statement. We include the proof, which is not difficult.

\begin{lemma}\label{lem:delta}
Let $f:\B^n\to\B^n$ and suppose $f$ has two distinct fixed points $x$ and $y$. Then $G(f)$ has a positive cycle $C$ such that $x_i\neq y_i$ for every vertex $i$ of $C$. 
\end{lemma}

\begin{proof}
Let $I$ be the set of components $i\in [n]$ with $x_i\neq y_i$, and, for all $i\in [n]$, let $s_i:=y_i-x_i$. Suppose that $s_i=1$. Seeking a contradiction, suppose that the following holds: for all $j\in [n]$, $x_j\geq y_j$ if $G(f)$ has a positive arc $j\to i$ and $x_j\leq y_j$ if $G(f)$ has a negative arc $j\to i$. It means that, going from $x$ to $y$, no activator of $i$ has been added, and no inhibitor of $i$ has been removed. It is then easy to check that this forces the inequality $f_i(x)\geq f_i(y)$. But since $s_i=1$, we have $f_i(x)=x_i<y_i=f_i(y)$ and we obtain the desired contradiction. This means that $G(f)$ has at least one positive arc $j\to i$ with $x_j<y_j$ or at least one negative arc $j\to i$ with $x_j>y_j$. In other words, $G(f)$ has an arc $j\to i$ of sign $s_js_i$. If $s_i=-1$, the same conclusion is obtained with similar arguments. Therefore, for every $i\in I$, there exists $j\in I$ such that $G(f)$ has an arc $j\to i$ of sign $s_js_i$. Hence, $G(f)$ has a cycle $C=i_0\to i_1\to\cdots\to i_{\ell-1}\to i_0$ where each vertex $i_k$ is in $I$, and where the sign of each arc $i_k\to i_{k+1}$ is $s_{i_k}s_{i_{k+1}}$, where $k+1$ is computed modulo~$\ell$. Hence, the sign of $C$ is (again by computing $k+1$ modulo $\ell$): 
\[
\prod_{k=0}^{\ell-1} s_{i_k}s_{i_{k+1}}
=\left(\prod_{k=0}^{\ell-1} s_{i_k}\right)\cdot \left(\prod_{k=0}^{\ell-1} s_{i_{k+1}}\right)
=\left(\prod_{k=0}^{\ell-1} s_{i_k}\right)\cdot \left(\prod_{k=0}^{\ell-1} s_{i_k}\right)
=1.
\]
Thus $C$ is positive and since all its vertices are in $I$ the lemma is proved. 
\end{proof}

In the spirit of Shih and Dong's theorem, Remy, Ruet and Thieffry \cite{RRT08} proved the following local version of Thomas' first rule. It can be regarded as the Boolean analogue of Cinquin-Demongeot's and Soul\'e's theorems \cite{CD02,S03} for differential equation systems. 

\begin{theorem}[\BF{Thomas' first rule, local version}]\label{thm:thomas1local}
Let $f:\B^n\to\B^n$ and suppose that $Gf(x)$ has no positive cycle for every $x\in\B^n$. Then $f$ has at most one fixed point.  
\end{theorem}

This generalizes Theorem~\ref{thm:thomas1} since if $G(f)$ has no positive cycle, then $Gf(x)$ has no positive cycle for every $x$, again since each $Gf(x)$ is a subgraph of $G(f)$. 

\smallskip
The following theorem draws a conclusion on $\Gamma(f)$. It shows that, in the absence of local positive cycles, there is a unique asynchronous attractors, reachable from every initial state by an  asynchronous geodesic. It remains valid in the discrete case \cite{RC07}. Thus positive cycles are not only necessary for multiple fixed points, but, more generally, for multiple asynchronous attractors.

\begin{theorem}
Let $f:\B^n\to\B^n$ and suppose that $Gf(x)$ has no positive cycle for every $x\in\B^n$. Then $\Gamma(f)$ has a unique attractor $A$. Furthermore, for every $x\in\B^n$, there exists $y\in A$ such that $\Gamma(f)$ has a geodesic from $x$ to $y$. 
\end{theorem}

\paragraph{Interlude} The above theorems show the importance of positive cycle. Then, is it difficult to check if a signed graph $G$ has a positive cycle? Let $G'$ be the signed graph obtained from $G$ by replacing each positive arc by a path of length two, with two negative arcs, and where the internal vertex is new. Let $D$ be the (unsigned) graph obtained from $G'$ by ignoring signs. If $G$ has $n$ vertices, then $D$ has at most $n^2+n$ vertices, and it is easy to check that $G$ has a positive cycle if and only if $D$ has an even cycle (a cycle of even length). The question thus reduces to the following: is it difficult to check if a graph $D$ has an even cycle? In other words, does there exist a polynomial time algorithm to check if $D$ has an even cycle? It turns out that this is a very difficult question, which remained open many years. Robertson, Seymour and Thomas \cite{RST99} finally proved the existence of a polynomial time algorithm for this decision problem. The algorithm designed is so elaborated that, up to our knowledge, no implementation exists. 

\bigskip
Let us say that a signed graph is \BF{negative} if it has no positive cycle. So when $G(f)$ is negative,  $f$ has at most one fixed point. Seeing that, it is natural to think about generalizations of the form {\em ``If $G(f)$ is not so far from being negative, then $f$ has not too much fixed points''}. We then need to formalize what we mean by {\em ``not so far from being negative''}, that is, we need a kind of distance from the negative case. For that, the notion of {\em positive transversal number}, defined below, is particularly suited. 

\smallskip
Let $G$ be a signed graph. A \BF{feedback vertex set} (FVS) of $G$ is a set of vertices $I$ intersecting every cycle (equivalently, removing the vertices in $I$ leaves the graph acyclic). The \BF{transversal number} $\tau$ of $G$ is the minimum size of a FVS of $G$. For instance, $\tau=0$ if and only if $G$ is acyclic. Similarly, a \BF{positive feedback vertex set} is a set of vertices $I$ intersecting every {\em positive} cycle (equivalently, removing the vertices in $I$ leaves the graph negative). The \BF{positive transversal number} $\tau^+$ of $G$ is the minimum size of a positive FVS of $G$. For instance, $\tau^+=0$ if and only if $G$ is negative. Note that $\tau^+\leq \tau$ since every FVS is a positive FVS.

\smallskip
The generalization we seek, proved by Aracena \cite{A08}, can now be stated as follows. The proof is a first rather easy application of Lemma~\ref{lem:delta}. 

\begin{theorem}[\BF{Positive feedback bound, global version}]
Let $f:\B^n\to\B^n$ and suppose that $I$ is a positive feedback vertex set of $G(f)$. Then $f$ has at most $2^{|I|}$
 fixed points. 
\end{theorem}

\begin{proof}
Let $X$ be the set of fixed points of $f$, and let $h$ be the map from $X$ to $\B^n$ defined by $h(x)_i=x_i$ for all $i\in I$ and $h(x)_i=0$ for all $i\in [n]\setminus I$. Thus the image of $h$ is of size at most $2^{|I|}$, that is, $|h(X)|\leq 2^{|I|}$. To completes the proof, it is thus sufficient to prove that $h$ is an injection. Let $x,y\in X$ be distinct fixed points. By Lemma~\ref{lem:delta}, $G(f)$ has a positive cycle $C$ such that $x_i\neq y_i$ for every vertex $i$ of $C$. By definition, $I$ intersects $C$. Hence, $C$ contains at least one vertex $i$ that belongs to $I$, and since $x_i\neq y_i$ we obtain $h(x)_i\neq h(y)_i$, so $h(x)\neq h(y)$. Thus $h$ is indeed an injection, and thus $|X|= |h(X)|\leq 2^{|I|}$, as desired. 
\end{proof}

Let $\tau$ and $\tau^+$ be the transversal and positive transversal numbers of $G(f)$. The positive feedback bound says that $f$ has at most $2^{\tau^+}$ fixed points. In particular, if $G(f)$ has no positive cycle, then $\tau^+=0$ and we recover Thomas' first rule: $f$ has at most one fixed point. Note that $2^\tau$ is a weaker upper-bound (since $\tau^+\leq\tau$) which uses less information on the interaction graph (signs are not taken into account). This weaker bound has been independently obtained by Riis \cite{R07} (see also \cite{GR11}), in the discrete setting, and in a completely different context: the network coding problem in Information Theory. 


\smallskip
The positive feedback bound remains valid in the non-Boolean discrete case and admits a local version that draws a conclusion not only on the fixed points, but more generally on the asynchronous attractors \cite{R09}. In the Boolean case, the statement is as follows. It implies both the positive feedback bound and the local version of Thomas' first rule.

\begin{theorem}[\BF{Positive feedback bound, local version}]
Let $f:\B^n\to\B^n$ and suppose that $I$ is a positive feedback vertex set of $Gf(x)$ for all $x\in\B^n$. Then $\Gamma(f)$ has at most $2^{|I|}$ attractors.
\end{theorem}

Let us finally discuss another application of Lemma~\ref{lem:delta}, given in \cite{GRR15} in the discrete setting. Fix an integer $d\geq 0$, and define $A(n,d)$ to be the maximum size of a subset $X$ of $\B^n$ such that $d(x,y)\geq d$ for all distinct $x,y\in X$. This quantity has been deeply studied in the context of Coding Theory. It is called the {\em maximum size of a binary code of length $n$ and minimum Hamming distance $d$}. In particular, the well known Gilbert's bound and sphere packing bound give the following approximation (see \cite{MWS77} for instance):
\[
\frac{2^n}{\sum_{k=0}^{d-1}{n \choose k}}\leq A(n,d)\leq \frac{2^n}{\sum_{k=0}^{\lfloor\frac{d-1}{2}\rfloor}{n \choose k}}. 
\]

\begin{theorem}[\BF{Coding bound}]
Let $f:\B^n\to\B^n$ and suppose that all the positive cycles of $G(f)$ are of length at least $d$. Then $f$ has at most $A(n,d)$ fixed points. 
\end{theorem}

\begin{proof}
Let $X$ be the set of fixed points of $f$. Let $x$ and $y$ be two distinct fixed points. Then, by Lemma~\ref{lem:delta}, $G(f)$ has a cycle $C$, say of length $\ell$, such that $x_i\neq y_i$ for every vertex $i$ of $G$. Thus $d(x,y)\geq\ell\geq d$. Consequently, $X$ is a binary code of length $n$ and minimum Hamming distance $d$, and thus $|X|\leq A(n,d)$. 
\end{proof}

The \BF{positive girth} of a signed graph is the minimum length of a positive cycle, and $\infty$ if no such cycle exists. The coding bound says that $f$ has at most $A(n,g^+)$ fixed points, where $g^+$ is the positive girth of $G(f)$. If $G(f)$ has no positive cycle, then $g^+=\infty$, and we deduce from the Gilbert's and sphere packing bounds that $A(n,g^+)=1$ (using the convention that ${n\choose k}=0$ for $k>n$). Thus $f$ has at most one fixed point, and we recover Thomas' first rule. See \cite{GR11,GRR15} for other connections with Coding Theory. 

\smallskip
At this stage, we have two upper-bounds on the fixed points of $f$: the positive feedback bound $2^{\tau^+}$ and the coding bound $A(n,g^+)$. It is thus natural to try to compare them. Clearly, they can be equal, as seen above if $G(f)$ has no positive cycle. But the positive feedback bound can be (and is often) much smaller than the coding bound. For instance, if $G$ is the disjoint union of a positive cycle of length one and a positive cycle of length $n-1$, then $f$ has $4$ fixed points. In that case $\tau^+=2$ and $g^+=1$. Hence, the positive feedback bound $2^{\tau^+}=4$ is tight while the coding bound $A(n,g^+)=2^n$ does not provide any information! We may think that the positive feedback bound is always at most the coding bound. We left this as an open question.

\begin{question}
Let $\tau^+$ and $g^+$ be the positive transversal number and positive girth of a signed graph $G$. Do we necessarily have $2^{\tau^+}\leq A(n,g^+)$?
\end{question}

\section{Thomas' second rule and beyond}

The following theorem, proved in \cite{R10} in the discrete case,  is a natural Boolean version of Thomas' second rule. It shows that negative cycles are necessary for asynchronous cyclic attractors. 

\begin{theorem}[\BF{Thomas' second rule, global version}]\label{thm:thomas2}
Let $f:\B^n\to\B^n$ and suppose that $G(f)$ has no negative cycle. Then all the attractors of $\Gamma(f)$ are fixed points. 
\end{theorem}

We may think that the conclusion ``all the attractors of $\Gamma(f)$ are fixed points'' is very strong. But it is (surprisingly) not so much. To see that let $\phi(n)$ be, among the $2^{n2^n}$ Boolean networks with $n$ components, the fraction of networks $f$ such that all the attractors of $\Gamma(f)$ are fixed points. Bollob\'as, Gotsman and Shamir \cite{BGS93} proved that $\phi(n)$ is, for large $n$, greater than~$\frac{1}{2}$. 

\begin{theorem} 
$\lim_{n\to\infty} \phi(n)=1-\frac{1}{e}$.
\end{theorem}

Since $\Gamma(f)$ obviously has at least one attractor, we obtain from Theorem~\ref{thm:thomas2} the following fixed point theorem. 

\begin{theorem}\label{thm:thomas2fixedpoint}
Let $f:\B^n\to\B^n$ and suppose that $G(f)$ has no negative cycle. Then $f$ has at least one fixed point.
\end{theorem}

Together, Theorems \ref{thm:thomas1} and \ref{thm:thomas2fixedpoint} give a nice ``proof by dichotomy'' of Robert's fixed point theorem, the first item in Theorem~\ref{thm:robert}. To see that, suppose that $G(f)$ is acyclic. Then $G(f)$ has no positive cycle, so $f$ has at most one fixed point (Theorem~\ref{thm:thomas1}), and $G(f)$ has no negative cycle, so $f$ has at least one fixed point (Theorem~\ref{thm:thomas2fixedpoint}). Thus $f$ has indeed a unique fixed point. 

\smallskip
Aracena \cite{A08} proved the following theorem which, compared with Theorem~\ref{thm:thomas2fixedpoint}, has a stronger conclusion under a stronger assumption (it is an easy exercise to deduce Theorem~\ref{thm:thomas2fixedpoint} from it, by considering the strongly connected components in the topological order). 

\begin{theorem}
Let $f:\B^n\to\B^n$ and suppose that $G(f)$ is strongly connected, has at least one arc, and has no negative cycle. Then $f$ has two fixed points $x$ and $y$ such that $d(x,y)=n$. 
\end{theorem}

We may now ask for a local version of Thomas' second rule, that is: {\em Is it true that if $Gf(x)$ has no negative cycle for all $x$, then $\Gamma(f)$ has no cyclic attractor?} A positive answer would provide, together with Theorem~\ref{thm:thomas1local}, a nice ``proof by dichotomy'' of Shih and Dong's theorem (Theorem~\ref{thm:shihdong}), in the spirit of the ``proof by dichotomy'' of Robert's theorem described above. In the non-Boolean discrete case, however, a negative answer was given in \cite{R10}, which leaded to a counter example of a local version of Thomas' second rule in the continuous case \cite{RC11}. After partial positive answers in the Boolean case \cite{R11,RR13}, Boolean counter examples were finally obtained independently, and almost simultaneously, by Ruet \cite{Ru2017} (with $n=7$), Tonello \cite{T17} (with $n=6$), and by Faur\'e and Kaji \cite{FK18} (with $n=6$). The two counter examples found with $n=6$ are different Boolean conversions of the non-Boolean discrete counter example given in \cite{R10}. Ruet's counter example is built with another technic. Recently, Tonello, Chaouiya and Farcot \cite{TFC18} proved, using a SAT solver, that there is no counter example with $n\leq 5$. 

\begin{theorem}[\BF{Counter example to the local version of Thomas' second rule}]
For all $n\geq 6$, there exists $f:\B^n\to\B^n$ such that $Gf(x)$ has no negative cycle for all $x\in\B^n$ and such that $f$ has no fixed point (this obviously implies that $\Gamma(f)$ has a cyclic attractor). 
\end{theorem}

We have seen that if $G(f)$ has no negative cycle, then all asynchronous attractors are fixed points. Can these fixed points be reached quickly? We have a very contrasted answer. Indeed, it is proven in \cite{MRRS13} that, for all states $x$ there exists a fixed point $y$ such that $\Gamma(f)$ has a geodesic from $x$ to $y$ (we recently discovered that Alon \cite{A85} obtained a similar result for threshold Boolean networks with positive edges weights). Thus, a quick convergence is always possible (this is Theorem~\ref{thm:MRRS13} below). However, it is proven in \cite{MRRS16} that, given an initial state $x$, it may exist a fixed point $y$, reachable from $x$ in $\Gamma(f)$, but with paths of exponential length only (this is Theorem~\ref{thm:MRRS16} below). Thus, some fixed points may be very hard to reach. 

\begin{theorem}\label{thm:MRRS13}
Let $f:\B^n\to\B^n$ and suppose that $G(f)$ has no negative cycle. Then, for every $x\in\B^n$, there exists a fixed point $y$ of $f$ such that $\Gamma(f)$ has a geodesic from $x$ to $y$. 
\end{theorem}

\begin{theorem}\label{thm:MRRS16}
For every $n\geq 1$, there exists $f:\B^n\to\B^n$ and $x,y\in\B^n$ with the following properties: $G(f)$ has only positive arcs, $y$ is a fixed point of $f$ reachable from $x$ in $\Gamma(f)$, and all the paths from $x$ to $y$ in $\Gamma(f)$ are of length at least $2^{\lfloor\frac{n}{2}\rfloor}$.
\end{theorem}

Note that Theorem~\ref{thm:MRRS13} trivially implies Theorem~\ref{thm:thomas2}. Indeed, suppose that $G(f)$ has no negative cycle, and suppose, for a contradiction, that $\Gamma(f)$ has a cyclic attractor $A$. Then no fixed point can be reached from $A$, and this contradicts Theorem~\ref{thm:MRRS13}. Thus all the attractors of $\Gamma(f)$ are fixed points, and we recover Theorem~\ref{thm:thomas2}. 

\paragraph{Interlude} Is it difficult to check if a signed graph $G$ has a negative cycle? We may assume, without loss, that $G$ is strongly connected (otherwise just consider, independently, its strongly connected components). Let $D$ be the graph obtained from $G$ as in the interlude of the previous section. Then $D$ is also strongly connected, and it is easy to check that $G$ has a negative cycle if and only if $D$ has an odd cycle (a cycle of odd length), which is equivalent to say that $D$ is not bipartite. The question thus reduces to the following: is it difficult to check if a strongly connected graph $D$ is bipartite?  There is an easy $O(n^2)$ time algorithm to check that. Indeed, it is sufficient to construct a spanning tree $T$ of $D$, using a depth-first search, and to consider a proper $2$-coloring of $T$. Then $D$ is bipartite if and only if this $2$-coloring is a proper coloring~of~$D$. 

\bigskip
What about the number of fixed points when there is no negative cycle? We have obviously the positive feedback-bound $2^{\tau^+}$, discussed in the previous section. Can this bound be improved in the absence of negative cycle? The answer is positive, and involves interesting tools from Set Theory. The key observation, given in \cite{MRRS13}, is the following. If $G$ is a signed graph, we denote by $|G|$ the underlying unsigned graph, obtained by ignoring signs. 

\begin{theorem}
Let $f:\B^n\to\B^n$ and suppose that $G(f)$ is strongly connected and has no negative cycle. Then there exists $h:\B^n\to\B^n$ with the following three properties: $G(h)$ has only positive arcs, $|G(h)|=|G(f)|$, and $\Gamma(h)$ is isomorphic to $\Gamma(f)$. 
\end{theorem}

\smallskip
This shows that, at least in the strongly connected case, the absence of negative cycle is equivalent to the absence of negative arc. This is a huge simplification, since if $G(f)$ has only positive arcs, then $f$ is very well structured. To see that, let $\leq$ be the partial order on $\B^n$ defined by $x\leq y$ if and only if $x_i\leq y_i$ for all $i\in [n]$. We say that $f:\B^n\to\B^n$ is \BF{monotone} if it preserves this partial order, that is, if $f(x)\leq f(y)$ for all states $x,y$ such that $x\leq y$. It is then easy to see that 
\[
\text{$f$ is monotone} \iff \text{$G(f)$ has only positive arcs.}
\]
Thus, at least in the strongly connected case, the absence of negative cycle is equivalent to the monotony. This leads us to focus on the monotone case. 

\smallskip
The first important theorem about monotone networks is the very well known Knaster-Tarski theorem \cite{KT28}, established in 1928. The statement needs a definition. A subset $L$ of $\B^n$ is a \BF{lattice} if, for all $x,y\in L$, the set of $z\in L$ such that $x\leq z$ and $y\leq z$ has a unique minimal element, and the set of $z\in L$ such that $z\leq x$ and $z\leq y$ has a unique maximal element. It is easy to check that every lattice has a unique maximal element and a unique minimal element. 

\begin{theorem}[\BF{Knaster-Tarski}]
If $f:\B^n\to\B^n$ is monotone, then the set of fixed points of $f$ is a non-empty lattice. 
\end{theorem}

This statement does not take into account the interaction graph of $f$, and improvements can be obtained by considering it. We again need some definitions, from Set and Graph Theory. Let $X\subseteq \B^n$ and $Y\subseteq \B^m$. We say that $X$ and $Y$ are \BF{isomorphic} if there exists a bijection $g:X\to Y$ such that, for all $x,y\in X$, $x\leq y$ if and only if $g(x)\leq g(y)$. A \BF{chain} in $X$ is a subset $C\subseteq X$ of pairwise comparable elements: for all $x,y\in C$, either $x\leq y$ or $y\leq x$. An \BF{antichain} in $X$ is a subset $A\subseteq X$ of pairwise incomparable elements: for all $x,y\in A$, neither $x\leq y$ nor $y\leq x$. Let $G$ be a signed graph. The \BF{packing number} $\nu$ of $G$ is the maximum size of a set of pairwise vertex-disjoint cycles. Similarly, the \BF{positive packing number} $\nu^+$ of $G$ is the maximum size of a set of pairwise vertex-disjoint {\em positive} cycles. Let $\tau$ and $\tau^+$ be the transversal and positive transversal numbers of $G$, respectively (defined in the previous section). It is easy to see that $\nu^+\leq\tau^+\leq \tau$ and $\nu^+\leq\nu\leq \tau$. Also, if $G$ has no negative cycle, then we obviously have $\tau^+=\tau$ and $\nu^+=\nu$. 

\smallskip
We are now in position to state an extension, proved in \cite{ARS17}, of the Knaster-Tarski theorem. 

\begin{theorem}\label{thm:monotone}
If $f:\B^n\to\B^n$ is monotone, then the set of fixed points of $f$ is isomorphic to a non-empty lattice $L\subseteq \B^\tau$ without chain of size $\nu+2$, where $\tau$ and $\nu$ be the transversal and packing numbers of $G(f)$, respectively.
\end{theorem}

Using the notations of the statement, and denoting by $X$ the set of fixed points of $f$, we have $|X|=|L|\leq 2^\tau$, thus we recover the positive feedback bound (since $G(f)$ has only positive arcs we have $\tau^+=\tau$). But the fact that $X$ is a lattice without chain of size $\nu+2$ provides informations on the structure of fixed points, which can be used to improve the bound. The improved bound is given in Theorem~\ref{thm:monotone_bound} below. The proof uses the following theorem of Erd\H{o}s \cite{E45} (the case $\ell=1$ is the well known Sperner's Lemma on the largest antichains of the $n$-cube \cite{S28}, namely: every antichain of $\B^n$ is of size at most ${n\choose \lfloor\frac{n}{2}\rfloor}$). 

\begin{theorem}[\BF{Erd\H{o}s}]
If $A$ is a subset of $\B^n$ without chain of size $\ell+1$, then $|A|$ is at most the sum of the $\ell$ largest binomial coefficients of order $n$, that is, 
\[
|A|\leq \sum_{k=\lfloor\frac{n-\ell+1}{2}\rfloor}^{\lfloor\frac{n+\ell-1}{2}\rfloor}{n\choose k}. 
\]
\end{theorem}

\begin{theorem}\label{thm:monotone_bound}
If $f:\B^n\to\B^n$ is monotone, then the number of fixed points of $f$ is at most two plus the the sum of the $\nu-1$ largest binomial coefficients of order $\tau$, where $\tau$ and $\nu$ be the transversal and packing numbers of $G(f)$, respectively. 
\end{theorem}

\begin{proof}
Using again the notations of Theorem~\ref{thm:monotone}, let $L$ be the subset $\B^\tau$ isomorphic to the set $X$ of fixed points of $f$. Since $L$ is a lattice, it has a unique maximal element, say $a$, and a unique minimal element, say $b$. Since all the chains of $L$ of maximum size contain both $a$ and $b$, and since $L$ has no chain of size $\nu+2$, we deduce that $L\setminus\{a,b\}$ has no chain of size $\nu$. Thus, according to Erd\H{o}s' theorem, $|L|-2$ is at most the sum of the $\nu-1$ largest binomial coefficients of order $\tau$. Since $|X|=|L|$, this proves the theorem.
\end{proof}

Let $\beta$ be the upper-bound given by the previous theorem. If $\nu=\tau$, then $\beta$ is two, {\em i.e.} ${\tau\choose 0}+{\tau\choose \tau}$, plus the sum of the $\tau-1$ remaining coefficients ${\tau \choose k}$. So $\beta$ is the sum of all the coefficients ${\tau\choose k}$, and thus $\beta=2^\tau$ is the feedback bound. But, when $\nu<\tau$, some coefficients are missing in the sum, and thus $\beta<2^{\tau}$ improves the feedback bound. Actually, the improvement is important when $\nu$ is much smaller than $\tau$. For instance, the complete graph on $n$ vertices (with an arc $j\to i$ if and only if $j\neq i$) has packing number $\nu=\lfloor\frac{n}{2}\rfloor$ and transversal number $\tau=n-1$, thus $\nu$ is linear with $\tau$. It is however very difficult to find graphs that exhibit larger gaps. The best construction is due to Alon and Seymour \cite{S95}, and shows that, for large $n$, we may have $\nu\log\nu\leq 30\tau$. It is also known that, for fixed $\nu$, $\tau$ cannot be arbitrarily large. This is a deep result of Reed, Robertson, Seymour and Thomas \cite{RRST96}, the following.

\begin{theorem}
There exists a function $r:\mathbb{N}\to\mathbb{N}$ such that, for every graph $G$, 
\[
\tau\leq r(\nu), 
\]
where $\tau$ and $\nu$ are the transversal and packing numbers of $G$, respectively. 
\end{theorem}

The function $r$ exhibited in the proof grows very fast, and the only exact value we know is $r(1)=3$ \cite{M93}. As a consequence, for every $f:\B^n\to\B^n$, the number of fixed points of $f$ is at most $2^{r(\nu)}$, where $\nu$ is the packing number of $G(f)$. Hence, the fixed points can be bound according to the maximum number of disjoint cycles only. Is it possible to bound the fixed points according to the maximum number of {\em positive} disjoint cycles only? This is a difficult open question, formally: 

\begin{question}
Does there exist a function $s:\mathbb{N}\to\mathbb{N}$ such that, for every $f:\B^n\to\B^n$, the number of fixed points of $f$ is at most $2^{s(\nu^+)}$, where $\nu^+$ is the positive packing number of~$G(f)$?
\end{question}

Let us finish this section with a third and last application of Lemma~\ref{lem:delta}, given in \cite{ARS17}, which partially proves Theorem~\ref{thm:monotone}, and which is connected to the above question. This shows again the usefulness of this lemma.

\begin{theorem}
If a function $f:\B^n\to\B^n$ has a chain of $\ell\geq 2$ fixed points, then $G(f)$ has at least $\ell-1$ vertex-disjoint positive cycles.
\end{theorem}

\begin{proof}
Suppose then that $f$ has a chain of $\ell\geq 2$ fixed points, say $x^1\leq x^2\leq\cdots\leq x^{\ell}$. For $1\leq k<\ell$, let $I^k$ be the set of components $i\in [n]$ such that $x^{k}_i\neq x^{k+1}_i$. By the chain property, the sets $I^k$ are pairwise distinct. Furthermore, by Lemma~\ref{lem:delta}, for each $1\leq k<\ell$, $G(f)$ has a positive cycle $C_k$ such that $x^k_i\neq x^{k+1}_i$ for every vertex $i$ in $C_k$. Thus the vertices of $C_k$ are all in $I^k$, and since the $I^k$ are pairwise disjoint, the $\ell-1$ positive cycles $C_k$ are pairwise disjoint. 
\end{proof}

\section{Connections between positive and negative cycles}

All the previous results consider negative cycles or positive cycles only. It's then rather clear that the next step consists in analyzing the two types of cycles simultaneously. This is difficult, and few results have been done in this direction. There is however, in \cite{R18}, the  following generalizations of Theorems~\ref{thm:thomas1} and \ref{thm:thomas2fixedpoint}, that take into account both positive and negative cycles. A strongly connected component is {\em non-trivial} if it contains at least one arc, and {\em initial} if there is no arc entering the component.  

\begin{theorem}\label{thm:R18}
Let $f:\B^n\to\B^n$ and suppose that every positive (resp. negative) cycle of $G(f)$ contains at least one arc whose deletion leaves a signed graph with a non-trivial initial strongly connected component containing only negative (resp. positive) cycles. Then $f$ has at most (resp. at least) one fixed point. 
\end{theorem} 

Suppose for instance that $G(f)$ is a kind of chain as follows:
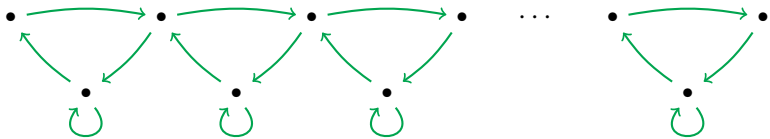
\[
\begin{tikzpicture}
\node (1) at (0,1){\scriptsize$\bullet$};
\node (2) at (2,1){\scriptsize$\bullet$};
\node (3) at (1,0){\scriptsize$\bullet$};
\node (4) at (4,1){\scriptsize$\bullet$};
\node (5) at (3,0){\scriptsize$\bullet$};
\node (6) at (6,1){\scriptsize$\bullet$};
\node (7) at (5,0){\scriptsize$\bullet$};
\node (...) at (7,1){$\cdots$};
\node (n3) at (8,1){\scriptsize$\bullet$};
\node (n1) at (10,1){\scriptsize$\bullet$};
\node (n) at (9,0){\scriptsize$\bullet$};
\draw[Green,->,thick] (3.-60) .. controls (1.5,-0.7) and (0.5,-0.7) .. (3.-120);
\draw[Green,->,thick] (5.-60) .. controls (3.5,-0.7) and (2.5,-0.7) .. (5.-120);
\draw[Green,->,thick] (7.-60) .. controls (5.5,-0.7) and (4.5,-0.7) .. (7.-120);
\draw[Green,->,thick] (n.-60) .. controls (9.5,-0.7) and (8.5,-0.7) .. (n.-120);
\path[Green,->,thick]
(1) edge[Green,bend left=10] (2)
(2) edge[Green,bend left=10] (3)
(3) edge[Green,bend left=10] (1)
(2) edge[Green,bend left=10] (4)
(4) edge[Green,bend left=10] (5)
(5) edge[Green,bend left=10] (2)
(4) edge[Green,bend left=10] (6)
(6) edge[Green,bend left=10] (7)
(7) edge[Green,bend left=10] (4)
(n3) edge[Green,bend left=10] (n1)
(n1) edge[Green,bend left=10] (n)
(n)  edge[Green,bend left=10] (n3)
;
\end{tikzpicture}
\]
There are $k=\frac{n-1}{2}$ loops (cycles of length one), and to remove all the positive cycles, it is necessary and sufficient to remove the $k$ vertices with a loop; for this graph: $\nu=\nu^+=\tau=\tau^+=k$. Hence, the positive feedback bound is $2^k$. This is rather good since it is easy to see that $f$ can have more than $2^{k-1}$ fixed points. Suppose now that $G(f)$ is as follows:
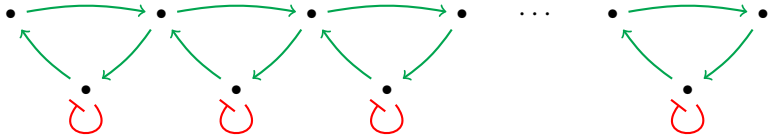
\[
\begin{tikzpicture}
\node (1) at (0,1){\scriptsize$\bullet$};
\node (2) at (2,1){\scriptsize$\bullet$};
\node (3) at (1,0){\scriptsize$\bullet$};
\node (4) at (4,1){\scriptsize$\bullet$};
\node (5) at (3,0){\scriptsize$\bullet$};
\node (6) at (6,1){\scriptsize$\bullet$};
\node (7) at (5,0){\scriptsize$\bullet$};
\node (...) at (7,1){$\cdots$};
\node (n3) at (8,1){\scriptsize$\bullet$};
\node (n1) at (10,1){\scriptsize$\bullet$};
\node (n) at (9,0){\scriptsize$\bullet$};
\draw[red,-|,thick] (3.-60) .. controls (1.5,-0.7) and (0.5,-0.7) .. (3.-120);
\draw[red,-|,thick] (5.-60) .. controls (3.5,-0.7) and (2.5,-0.7) .. (5.-120);
\draw[red,-|,thick] (7.-60) .. controls (5.5,-0.7) and (4.5,-0.7) .. (7.-120);
\draw[red,-|,thick] (n.-60) .. controls (9.5,-0.7) and (8.5,-0.7) .. (n.-120);
\path[red,->,thick]
(1) edge[Green,bend left=10] (2)
(2) edge[Green,bend left=10] (3)
(3) edge[Green,bend left=10] (1)
(2) edge[Green,bend left=10] (4)
(4) edge[Green,bend left=10] (5)
(5) edge[Green,bend left=10] (2)
(4) edge[Green,bend left=10] (6)
(6) edge[Green,bend left=10] (7)
(7) edge[Green,bend left=10] (4)
(n3) edge[Green,bend left=10] (n1)
(n1) edge[Green,bend left=10] (n)
(n)  edge[Green,bend left=10] (n3)
;
\end{tikzpicture}
\]
To remove all the positive cycles, it is necessary and sufficient to remove $k-1$ vertices; for this graph: $\tau^+=k-1$ (and $\nu^+=\lfloor\frac{k}{2}\rfloor$). Thus the positive feedback bound is $2^{k-1}$. However, according to Theorem~\ref{thm:R18}, $f$ has at most one fixed point. Indeed, take any positive cycle (necessarily of length three)  and remove the arc whose head has a negative loop: this negative loop becomes a non-trivial initial strongly connected component with only negative cycles, thus the hypothesis of Theorem~\ref{thm:R18} is satisfied. In that case, the feedback bound is very bad... and it is clear that the problem comes from the fact that negative cycles are not taken into account. 

\smallskip
Seeing the above two examples, we may think that negative cycles are harmful for the presence of many fixed points. But this is far from being true. To see that, let $K^+_n$ (resp. $K^-_n$) be the {\em full-positive} (resp. {\em full-negative}) {\em complete graph}: the vertex set is $[n]$ and there is a positive (resp. negative) arc from $i$ to $j$ for all distinct $i,j\in [n]$. Given a signed graph $G$, let $\max(G)$ be the maximum number of fixed points in a Boolean network $f$ with $G(f)=G$. The following, proven in \cite{GRR15}, is a rather immediate application of Sperner lemma (already mentioned) and two classical results from Coding Theory. 

\begin{theorem}
For all n$\geq 1$, 
\[
\frac{{n\choose \lfloor\frac{n}{2}\rfloor}}{n}
\leq
\max(K^+_n)
\leq 
\frac{2^{n+1}}{n+2} 
\leq
\max(K^-_n)
=
{n\choose \lfloor\frac{n}{2}\rfloor} \sim \frac{2^n}{\sqrt{\pi n/2}}.
\]
\end{theorem}

As a consequence, for any fixed $k$, if $n$ is sufficiently large, then there exists a network with $K^-_n$ as interaction graph which has {\em more} fixed points than {\em any} network with $K^+_{n+k}$ as interaction graph. However, all is greater in  $K^+_{n+k}$ concerning the positive cycles. In particular the positive feedback number is $\lfloor\frac{n+k}{2}\rfloor$ while those of $K^-_n$ is $\lfloor\frac{n}{2}\rfloor$. These are therefore the negative cycles of $K^-_n$ that allow the presence of many fixed points. See \cite{ARS14,ARS17b} for other results that show that negative cycles are necessary for the presence of many fixed points. All these results suggest that, in order to have many fixed points, we need something like many short positive cycles connected each other through long negative cycles. 

\smallskip
In summary, the influence of connections between positive and negative cycles seems very subtile and is a widely open subject.  

\section{Conclusion}

Thomas' rules and Thomas' logical description of gene networks leaded to a body of results, in discrete  mathematics, about the relationships between the interaction graph and the dynamical properties of a Boolean (or discrete) network. These results involve many techniques from other fields, such as Graph Theory, Coding Theory and Set Theory (and results obtained in the context of networks can feed in turn these fields, see for instance  \cite{GRF16,ARS17b}). 

\smallskip
We have focussed on the number of fixed points, about which there are several results, and also on the attractors and geodesic paths of the asynchronous graph. This should appear as the beginning of something bigger since many other dynamical properties should be considered, such as the number and length of limit cycles of the synchronous graph. There are several results in this direction \cite{GO81,PS83,RT85,G15rank}, but the influence of the signed interaction graph on these kind of  dynamical parameters is still widely unknown. We have also taken an ``extremal'' point of view, studying for instance the maximum number of fixed points, but average analysis should be also relevant, see for instance \cite{G15rank,G17} . 

\smallskip
It is also important to compare the different dynamics that can be derived from $f$, such as the synchronous and asynchronous one. Indeed, Kauffman \cite{K69,K93} also introduced Boolean networks to model gene networks, but by considering the synchronous dynamics. Since then, to model the same biological object, some researchers use the synchronous dynamics (which is not consistent with continuous models) while others use the asynchronous dynamics (which is consistent with continuous models). However, the two dynamics are generally very different (for instance, in the synchronous case, Thomas' second rule is false, as showed in Fig.~\ref{fig:positive_cycle}). It is therefore important to compare them. Concerning this, interesting results have been obtained in \cite{NS17}. 

\smallskip
There are actually many other directions that needed to be explored to hopeful converge toward a consistent theory of Boolean networks, and more generally of finite dynamical systems. Thomas' ideas will live longtime in this context.  
 
\paragraph{Acknowledgment}
I wish to thank Ren\'e Thomas for his encouragements and his kindness. I remember that after the talk I gave for his 80th birthday, he said to me ``I did not fully understood all the results you presented, but I understood they are all important''. This had a strong positive impact on me (even if I thought my presentation might not be so clear).

\bibliographystyle{plain}
\bibliography{BIB}

\end{document}